\documentclass[12pt]{amsart}
\usepackage{amssymb,amsmath,amsthm,latexsym,cancel}
\usepackage{graphicx}
\usepackage{epsfig}
\newtheorem{theorem}{Theorem}[section]
\newtheorem{lemma}[theorem]{Lemma}

\theoremstyle{remark}
\newtheorem{definition}[theorem]{Definition}
\newtheorem{remark}[theorem]{Remark}

\newcommand\bs{\backslash}

\begin{document}

\title{On some ternary operations in knot theory}
\date{November 6, 2012}
\author{Maciej Niebrzydowski}
\address[Maciej Niebrzydowski]{Department of Mathematics\\
	 University of Louisiana at Lafayette\\
	 1403 Johnston Street\\
	 217 Maxim D. Doucet Hall\\
	 Lafayette, LA 70504-1010}
\email{mniebrz@gmail.com}


\thispagestyle{empty}

\begin{abstract}
We introduce a way to color the regions of a classical knot diagram using ternary operations, so that the number of colorings is a knot invariant. By choosing appropriate substitutions in the algebras that we assign to diagrams, one obtains the relations from the knot group, and from the core group.
Using the ternary operator approach, we generalize the Dehn presentation of the knot group to extra loops, and a similar presentation for the core group to the variety of Moufang loops.
\end{abstract}

\maketitle

\section{Introduction}
Encouraged by the existence of the Dehn presentation of the knot group, we introduce a way to color the regions of a classical knot diagram using ternary operations, so that the number of colorings is a knot invariant. To every knot diagram, oriented or unoriented, we assign an algebra that we call the ternary algebra of a knot. It is an invariant under Reidemeister moves, and every coloring using ternary operations can be viewed as a homomorphism from the ternary algebra of the knot to the algebra used for the coloring.

We consider the case of unoriented diagrams first. We introduce the ternary algebras using checkerboard colorings, and explain how the relations in the knot group and in the core group can be described via ternary operations. Then we move on to ternary invariants for oriented diagrams. The penultimate section contains the definitions and the proof of invariance of the general ternary algebras associated to diagrams. In the last section, we search for ternary operators that use binary operations from some nonassociative structures.

\begin{figure}
\begin{center}
\includegraphics[height=1.2 cm]{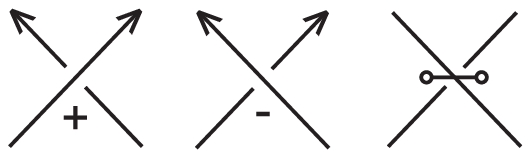}
\caption{}\label{conv}
\end{center}
\end{figure}

\begin{figure}
\begin{center}
\includegraphics[height=2 cm]{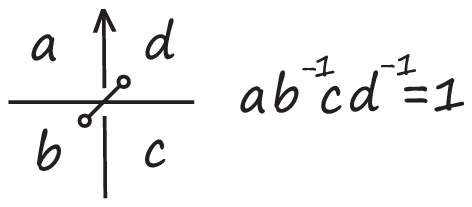}
\caption{}\label{Dehn}
\end{center}
\end{figure}

\section{Unoriented case}\label{unoriented}
We use the standard conventions from knot theory specifying the sign of a crossing, and a positive marker; see Figure \ref{conv}.
First, we recall, after \cite{Kau}, that the fundamental group of the complement of a knot in $\mathbb{S}^3$ can be given the following presentation, called the Dehn presentation: generators are assigned to the regions of a diagram, and relations correspond to crossings and are as in Figure \ref{Dehn}. One of the generators, say the one corresponding to the unbounded region, is set equal to identity. This definition uses the orientation of the under-arc at a crossing, but we note that it can be discarded if we use the positive marker instead. Consider a ternary operation $xyzT=xy^{-1}z$. A group with this operation is often given as an example of an algebra satisfying Mal'cev identities: $xyyT=x=yyxT$. The relation at a crossing is equivalent to any one of the following conditions: $a=dcbT$, $b=cdaT$, $c=badT$, and $d=abcT$. Thus, we see that each generator around a crossing can be expressed using the operator $T$ and the remaining generators, under the condition that the arguments for $T$ are read counter-clockwise if the region contains the positive marker, and clockwise otherwise.

\begin{figure}
\begin{center}
\includegraphics[height=2.6 cm]{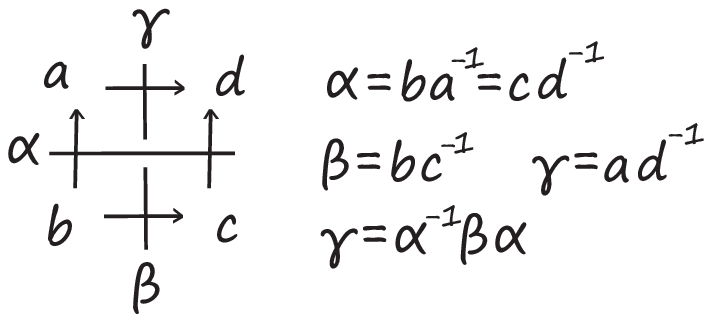}
\caption{}\label{Wirtinger}
\end{center}
\end{figure}

To go back to the standard Wirtinger relations of $\pi_1(\mathbb{S}^3\setminus L)$, we can proceed as follows. Given an unoriented diagram with generators satisfying the above relations, assign a co-orientation to the components of the link $L$; any co-orientation will do as long as it is consistent on each component. A standard generator $x$ assigned to an arc is obtained by setting $x=uv^{-1}$, where $u$ and $v$ are the regions adjacent to the arc labeled $x$, and the co-orientation points from $u$ to $v$. Thus, for the co-orientation arrows in Fig. \ref{Wirtinger}, we set $\alpha=ba^{-1}=cd^{-1}$, $\beta=bc^{-1}$, and $\gamma=ad^{-1}$. Then $\alpha^{-1}\beta\alpha=(ba^{-1})^{-1}(bc^{-1})(cd^{-1})=ad^{-1}=\gamma,$ as needed.

\begin{figure}
\begin{center}
\includegraphics[height=2.5 cm]{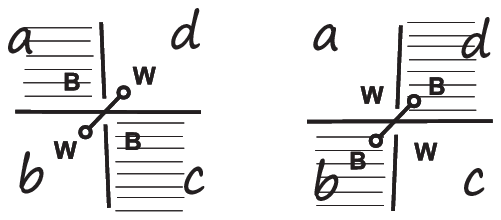}
\caption{}\label{unorientedop}
\end{center}
\end{figure}

Now we will generalize this idea. Let $D$ be a checkerboard-colored diagram of a link. We fix the convention that the unbounded region is white. Let $\mathcal{A}$ be an algebra with two ternary operators: $B$ for the black and $W$ for the white regions.

\begin{definition}\label{colunor}
A {\it checkerboard ternary coloring} is an assignment of the elements of an algebra $\mathcal{A}$ to the regions of the diagram in such a way that, at a crossing, if a generator corresponds to a white (resp. black) region then it is expressed in terms of the other three generators using the operator $W$ (resp. $B$). If the region contains the positive marker, then the inputs for the operators are read counter-clockwise, otherwise they are taken clockwise. More precisely, in the left part of Fig. \ref{unorientedop}, we have the following relations: $a=dcbB$, $b=cdaW$, $c=badB$, and $d=abcW$. In the situation on the right of Fig. \ref{unorientedop}, the equations are: $a=dcbW$, $b=cdaB$, $c=badW$, and $d=abcB$.
\end{definition}

We will now check what conditions have to be satisfied by the algebra $\mathcal{A}$ so that the number of the above colorings becomes a link invariant.
For fixed elements $a$, $b\in\mathcal{A}$, we consider several maps from $\mathcal{A}$ to $\mathcal{A}$: $B_{a,x,b}=axbB$, $B_{x,a,b}=xabB$,
$B_{a,b,x}=abxB$, $W_{a,x,b}=axbW$, $W_{x,a,b}=xabW$, and $W_{a,b,x}=abxW$. By making substitutions in the equations in the Definition \ref{colunor}, for example $a=dcbB=(abcW)cbB$, we obtain the following conditions:
\begin{enumerate}
\item[(1)] $B_{a,x,b}$ is inverse to $B_{b,x,a}$ ($b(axbB)aB=x$ and $a(bxaB)bB=x$)
\item[(2)] $W_{a,x,b}$ is inverse to $W_{b,x,a}$
\item[(3)] $B_{a,b,x}$ is inverse to $W_{b,a,x}$
\item[(4)] $B_{x,a,b}$ is inverse to $W_{x,b,a}$.
\end{enumerate}

\begin{figure}
\begin{center}
\includegraphics[height=3 cm]{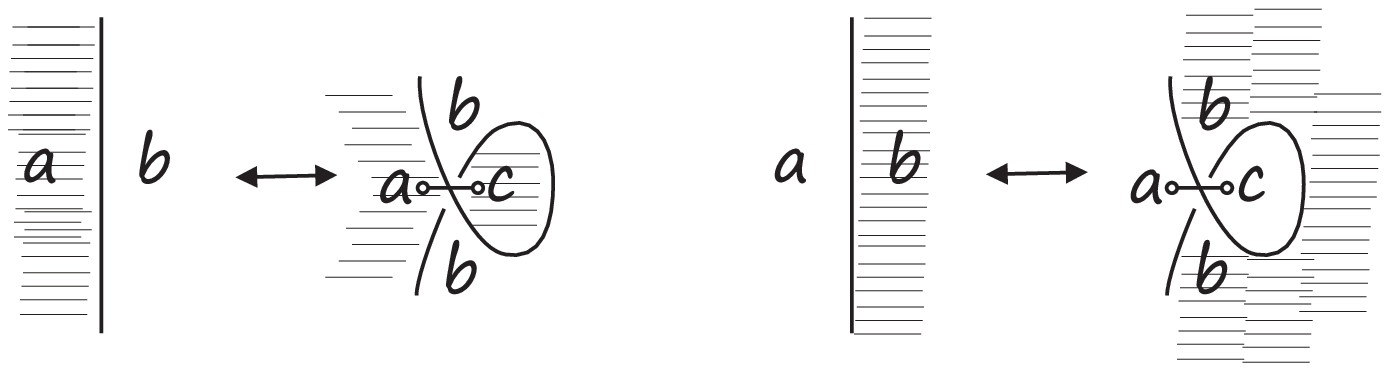}
\caption{}\label{R1unor}
\end{center}
\end{figure}

We note that for the first Reidemeister move (see Fig. \ref{R1unor}), the number of colorings does not change, because the color of the new region is determined by the operators: $c=babB$ on the left, and $c=babW$ on the right of the Fig. \ref{R1unor}.

\begin{figure}
\begin{center}
\includegraphics[height=4.3 cm]{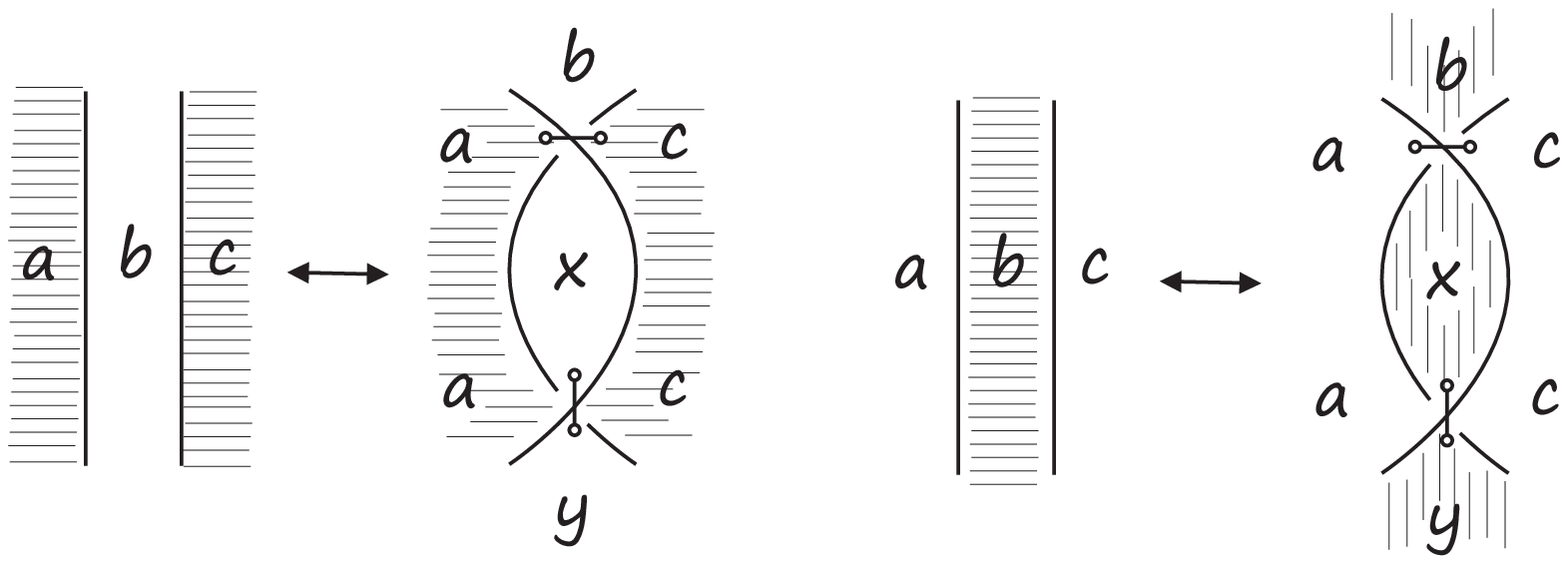}
\caption{}\label{R2unor}
\end{center}
\end{figure}

In the case of the second Reidemeister moves (Fig. \ref{R2unor}), the invariance follows from the conditions (1) and (2): on the left we have $x=abcW$, and $y=cxaW=c(abcW)aW=b$; on the right $x=abcB$ and $y=cxaB=c(abcB)aB=b$.

\begin{figure}
\begin{center}
\includegraphics[height=4.5 cm]{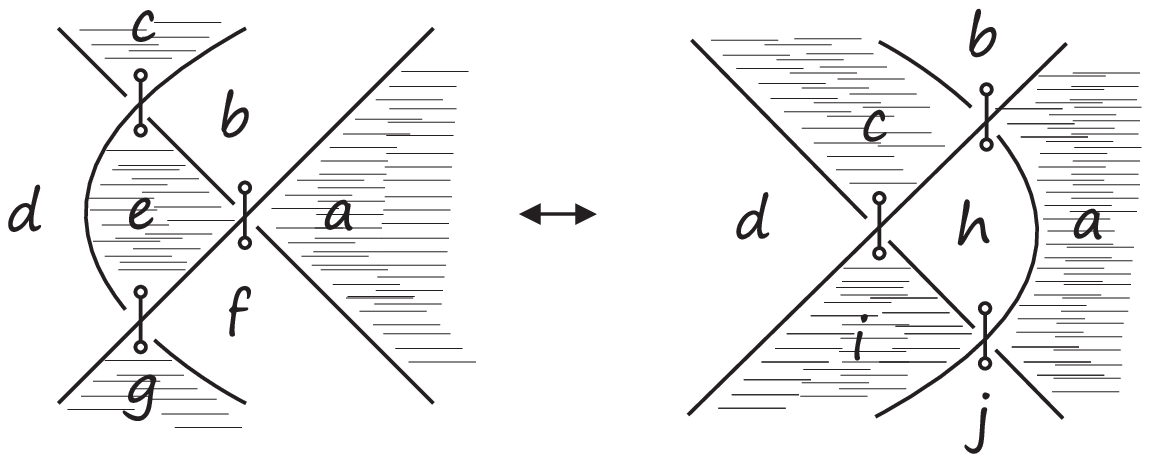}
\caption{}\label{R3unor}
\end{center}
\end{figure}

From the third Reidemeister move, we get relations that perhaps are the three-variable counterparts of self-distributivity for binary operations. With labeling as in Fig. \ref{R3unor}, on the left we have: $e=bcdB$, $f=abeW=ab(bcdB)W$, $g=fedB=[ab(bcdB)W](bcdB)dB$; on the right:
$h=abcW$, $i=hcdB=(abcW)cdB$, $j=ahiW=a(abcW)[(abcW)cdB]W$. Because the color $g$ should be equal to $i$, and $f$ equal to $j$, we get the following relations: 
\begin{enumerate}
\item[(5)] $(abcW)cdB=[ab(bcdB)W](bcdB)dB$; note that the right side of this equation is obtained from the left side by substituting $bcdB$ for $c$.
\item[(6)] $ab(bcdB)W=a(abcW)[(abcW)cdB]W$. Here, replacing the doubled letter $b$ by $abcW$ gives the right hand side.
\end{enumerate}
The possibility of the opposite checkerboard coloring, makes it necessary to add two more axioms with $B$ replaced by $W$ and vice versa:
\begin{enumerate}
\item[(7)] $(abcB)cdW=[ab(bcdW)B](bcdW)dW$
\item[(8)] $ab(bcdW)B=a(abcB)[(abcB)cdW]B$.
\end{enumerate}

It follows from the axioms (1)-(4) that the operators $W$ and $B$ can be presented (for finite algebras) as Latin cubes.
An example of a four-element algebra $\mathcal{A}$ satisfying the conditions (1) through (8) is given below ($x$ corresponds to rows and $y$ to columns). 

\begin{equation*}
xy4W=
\begin{bmatrix}
1 & 2 & 3 & 4 \\
3 & 4 & 2 & 1 \\
4 & 3 & 1 & 2 \\
2 & 1 & 4 & 3 
\end{bmatrix} \qquad
xy4B=
\begin{bmatrix}
3 & 1 & 4 & 2 \\
2 & 4 & 1 & 3 \\
1 & 2 & 3 & 4 \\
4 & 3 & 2 & 1 
\end{bmatrix} 
\end{equation*}
\begin{equation*}
xy3W=
\begin{bmatrix}
4 & 3 & 1 & 2 \\
1 & 2 & 3 & 4 \\
2 & 1 & 4 & 3 \\
3 & 4 & 2 & 1 
\end{bmatrix} \qquad
xy3B=
\begin{bmatrix}
2 & 4 & 1 & 3 \\
3 & 1 & 4 & 2 \\
4 & 3 & 2 & 1 \\
1 & 2 & 3 & 4 
\end{bmatrix} 
\end{equation*}

\begin{equation*}
xy2W=
\begin{bmatrix}
3 & 4 & 2 & 1 \\
2 & 1 & 4 & 3 \\
1 & 2 & 3 & 4 \\
4 & 3 & 1 & 2 
\end{bmatrix} \qquad
xy2B=
\begin{bmatrix}
1 & 2 & 3 & 4 \\
4 & 3 & 2 & 1 \\
2 & 4 & 1 & 3 \\
3 & 1 & 4 & 2 
\end{bmatrix}
\end{equation*} 
\begin{equation*}
xy1W=
\begin{bmatrix}
2 & 1 & 4 & 3 \\
4 & 3 & 1 & 2 \\
3 & 4 & 2 & 1 \\
1 & 2 & 3 & 4 
\end{bmatrix} \qquad
xy1B=
\begin{bmatrix}
4 & 3 & 2 & 1 \\
1 & 2 & 3 & 4 \\
3 & 1 & 4 & 2 \\
2 & 4 & 1 & 3 
\end{bmatrix} 
\end{equation*}

Now we consider the core group of a link \cite{FR,Joy}, $core(L)$. Its generators correspond to arcs of a diagram, and relations come from crossings, and are of the form $\gamma=\alpha\beta^{-1}\alpha$, where $\gamma$ and $\beta$ are assigned to the under-arcs of a crossing, and $\alpha$ is for the bridge. From the point of view of algebraic topology, $core(L)$ is the free product of the fundamental group of the cyclic branched double cover of $\mathbb{S}^3$ with branching set $L$ and the infinite cyclic group \cite{Wad, Prz}.

\begin{figure}
\begin{center}
\includegraphics[height=2.6 cm]{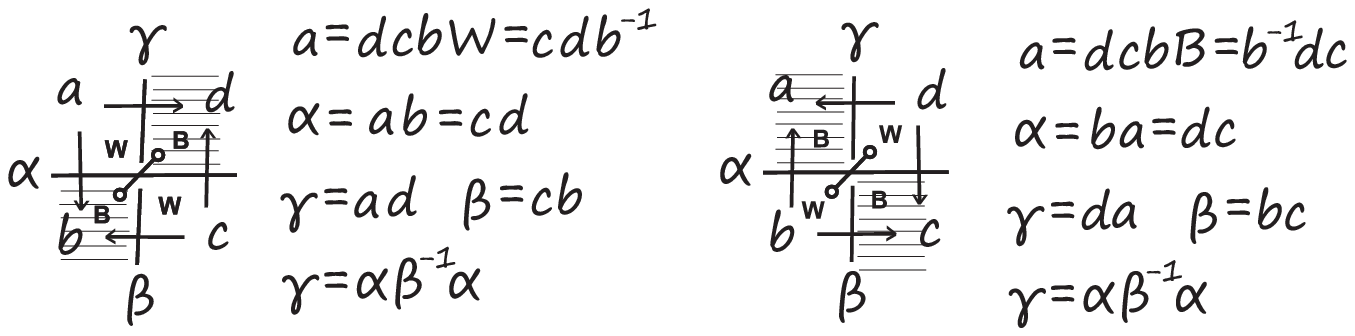}
\caption{}\label{core}
\end{center}
\end{figure}

Let $\mathcal{A}$ be a group, and the operators $B$ and $W$ be defined by $abcB=c^{-1}ab$, and $abcW=bac^{-1}$. The algebra $(\mathcal{A},B,W)$ satisfies the conditions (1) through (8). As in the case of the knot group, we can obtain the relations of the core group from the checkerboard ternary coloring. To each arc of the diagram we assign an element $x=uv$, where $u$ and $v$ are the elements assigned to, respectively, white and black regions adjacent to the arc labeled $x$. Then, in the situation as on the left of Fig. \ref{core}, we set $\alpha=ab=cd$ (note that $a=dcbW=cdb^{-1}$), $\beta=cb$, and $\gamma=ad$. It follows that $\alpha\beta^{-1}\alpha=(ab)(cb)^{-1}(cd)=ad=\gamma,$ as needed. For the coloring as on the right, we have: $\alpha=ba=dc$ (here, 
$a=dcbB=b^{-1}dc$), $\beta=bc$, and $\gamma=da$. It follows that $\alpha\beta^{-1}\alpha=(dc)(bc)^{-1}(ba)=da=\gamma,$ as desired.

A computer search (using GAP \cite{GAP4}) for group words giving ternary operators $B$ and $W$ yielded the following examples (the roles of $B$ and $W$ could be exchanged):
\begin{enumerate}
\item[(g1)] $B(a,b,c)=ab^{-1}c$, $W(a,b,c)=ab^{-1}c$; 
\item[(g2)] $B(a,b,c)=ac^{-1}b$, $W(a,b,c)=ac^{-1}b$;
\item[(g3)] $B(a,b,c)=ba^{-1}c$, $W(a,b,c)=ba^{-1}c$;
\item[(g4)] $B(a,b,c)=bc^{-1}a$, $W(a,b,c)=bc^{-1}a$;  
\item[(g5)] $B(a,b,c)=ca^{-1}b$, $W(a,b,c)=ca^{-1}b$;
\item[(g6)] $B(a,b,c)=cb^{-1}a$, $W(a,b,c)=cb^{-1}a$; 
\item[(g7)] $B(a,b,c)=a^{-1}cb$, $W(a,b,c)=bca^{-1}$;
\item[(g8)] $B(a,b,c)=c^{-1}ab$, $W(a,b,c)=bac^{-1}$; 
\item[(g9)] $B(a,b,c)=c^{-1}b^{-1}a^{-1}$, $W(a,b,c)=a^{-1}b^{-1}c^{-1}$.  
\end{enumerate}
Among the above pairs of operators, (g1), (g3), (g4), (g5), (g6), (g8) give the relations of the core group, and (g2), (g7), (g9) give the relations of the knot group. This observation will be useful when we consider the invariants obtained from loops.

\section{Oriented case}

\begin{figure}
\begin{center}
\includegraphics[height=2.8 cm]{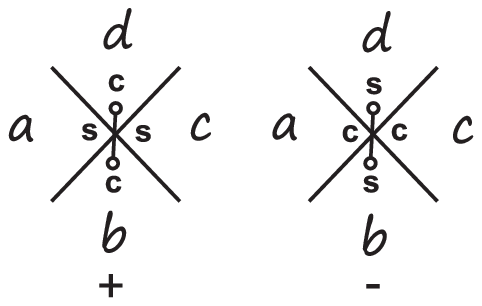}
\caption{}\label{orientedop}
\end{center}
\end{figure}

Now we define similar invariants for oriented links. 

\begin{definition}\label{col_or}
Let $\mathcal{A}$ be an algebra with ternary operators $S$ and $C$. A {\it ternary coloring} of a link diagram $D$ is an assignment of the elements of $\mathcal{A}$ to the regions of the diagram in such a way that at a positive crossing depicted on the left side of Fig. \ref{orientedop}, the relations are as follows:
$a=dcbS$, $b=cdaC$, $c=badS$, and $d=abcC$. Thus, for a positive crossing, the regions with a positive marker are assigned the operator $C$, and the other two regions get $S$; this is reversed if the crossing is negative. As in the unoriented case, the inputs for the operators are collected counter-clockwise if the region contains a positive marker, and clockwise otherwise. Thus, for the negative crossing on the right of Fig. \ref{orientedop}, we require the following relations: $a=dcbC$, $b=cdaS$, $c=badC$, and $d=abcS$.
\end{definition}

As before, appropriate substitutions give the set of axioms:

\begin{enumerate}
\item[(1)] $C_{a,x,b}$ is inverse to $C_{b,x,a}$
\item[(2)] $S_{a,x,b}$ is inverse to $S_{b,x,a}$
\item[(3)] $C_{a,b,x}$ is inverse to $S_{b,a,x}$
\item[(4)] $C_{x,a,b}$ is inverse to $S_{x,b,a}$.
\end{enumerate}

To check which conditions are imposed by the oriented Reidemeister moves, we use one of the generating sets of such moves found in \cite{Pol}. 

\begin{figure}
\begin{center}
\includegraphics[height=5 cm]{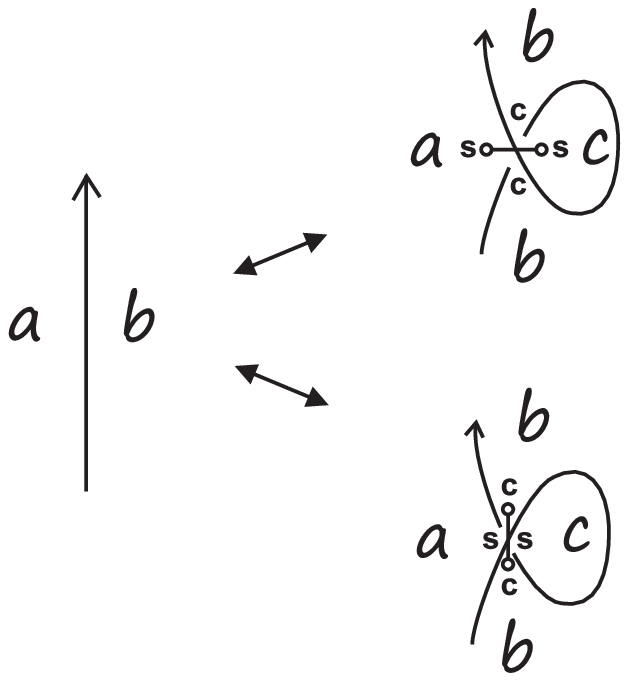}
\caption{}\label{R1or}
\end{center}
\end{figure}

\begin{figure}
\begin{center}
\includegraphics[height=4.3 cm]{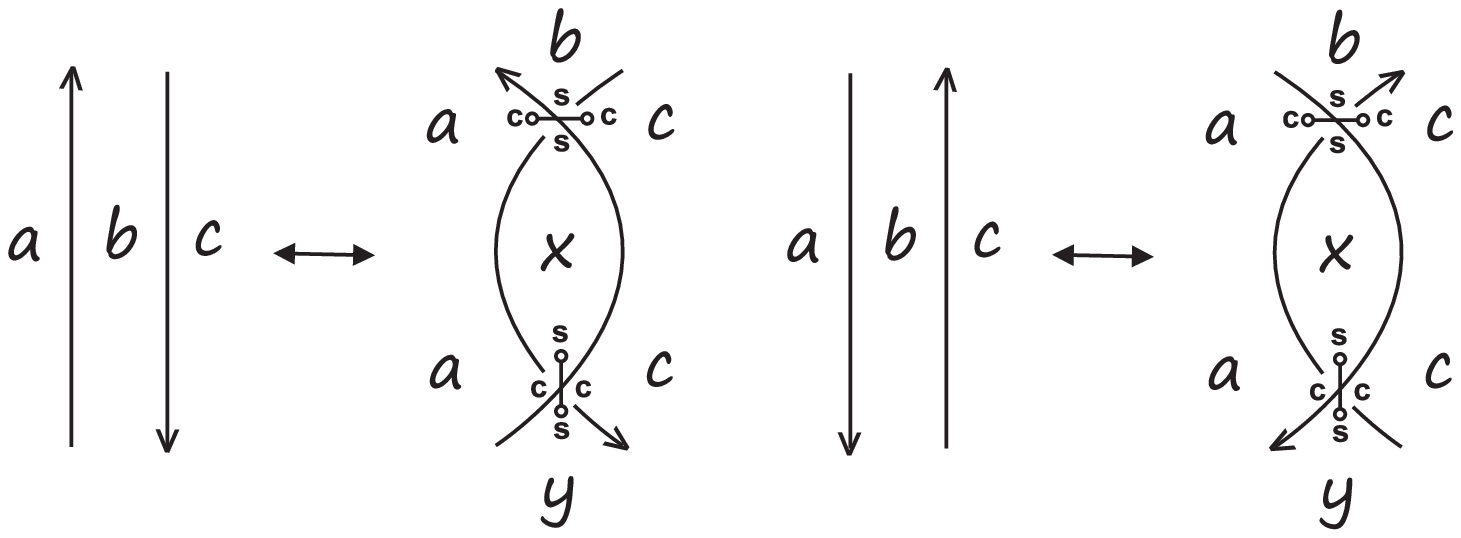}
\caption{}\label{R2or}
\end{center}
\end{figure}

The first Reidemeister move does not change the number of ternary colorings; in both versions of the first move, depicted in Fig. \ref{R1or}, $c=babC$. 

It turns out that for both second Reidemeister moves from the generating set, the calculation is the same (see Fig. \ref{R2or}): $x=abcS$, 
$y=cxaS=c(abcS)aS=b$. Thus, no new conditions are needed.

\begin{figure}
\begin{center}
\includegraphics[height=5 cm]{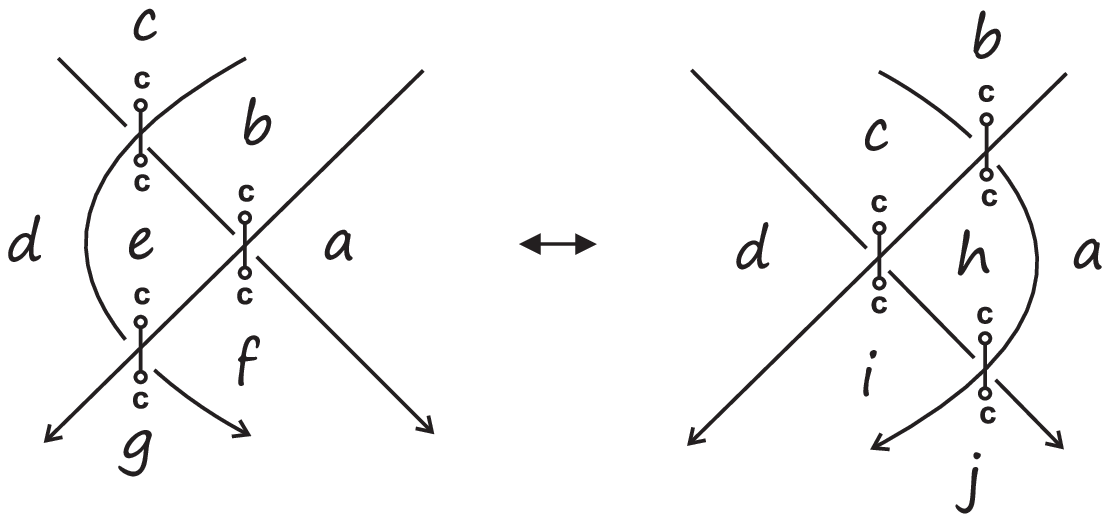}
\caption{}\label{R3or}
\end{center}
\end{figure}

The third Reidemeister move gives the `distributivity' axiom, this time involving only one operator. On the left of Fig. \ref{R3or} we have:
$e=bcdC$, $f=abeC=ab(bcdC)C$, $g=fedC=[ab(bcdC)C](bcdC)dC$. After the move: $h=abcC$, $i=hcdC=(abcC)cdC$, $j=ahiC=a(abcC)[(abcC)cdC]C$. By comparing $g$ with $i$, and $f$ with $j$, we get:
\begin{enumerate}
\item[(5)] $(abcC)cdC=[ab(bcdC)C](bcdC)dC$
\item[(6)] $ab(bcdC)C=a(abcC)[(abcC)cdC]C$.
\end{enumerate}

The computer search for group words giving operators $S$ and $C$ did not give any examples with $S\neq C$ ($S=C$ forgets about the orientation). We found the example $C(a,b,c)=S(a,b,c)=ab^{-1}c$, and all its relatives obtained by permuting the letters $a$, $b$, and $c$. It is however easy to generate finite examples with $S\neq C$. One such algebra is given below. Again, both operators can be represented by Latin cubes.

\begin{equation*}
xy4C=
\begin{bmatrix}
2 & 3 & 4 & 1 \\
1 & 4 & 3 & 2 \\
4 & 1 & 2 & 3 \\
3 & 2 & 1 & 4 
\end{bmatrix} \qquad
xy4S=
\begin{bmatrix}
2 & 3 & 4 & 1 \\
3 & 4 & 1 & 2 \\
4 & 1 & 2 & 3 \\
1 & 2 & 3 & 4 
\end{bmatrix} 
\end{equation*}
\begin{equation*}
xy3C=
\begin{bmatrix}
3 & 4 & 1 & 2 \\
4 & 1 & 2 & 3 \\
1 & 2 & 3 & 4 \\
2 & 3 & 4 & 1 
\end{bmatrix} \qquad
xy3S=
\begin{bmatrix}
3 & 2 & 1 & 4 \\
4 & 1 & 2 & 3 \\
1 & 4 & 3 & 2 \\
2 & 3 & 4 & 1 
\end{bmatrix} 
\end{equation*}

\begin{equation*}
xy2C=
\begin{bmatrix}
4 & 1 & 2 & 3 \\
3 & 2 & 1 & 4 \\
2 & 3 & 4 & 1 \\
1 & 4 & 3 & 2 
\end{bmatrix} \qquad
xy2S=
\begin{bmatrix}
4 & 1 & 2 & 3 \\
1 & 2 & 3 & 4 \\
2 & 3 & 4 & 1 \\
3 & 4 & 1 & 2
\end{bmatrix}
\end{equation*} 
\begin{equation*}
xy1C=
\begin{bmatrix}
1 & 2 & 3 & 4 \\
2 & 3 & 4 & 1 \\
3 & 4 & 1 & 2 \\
4 & 1 & 2 & 3  
\end{bmatrix} \qquad
xy1S=
\begin{bmatrix}
1 & 4 & 3 & 2 \\
2 & 3 & 4 & 1 \\
3 & 2 & 1 & 4 \\
4 & 1 & 2 & 3 
\end{bmatrix} 
\end{equation*}

\section{General ternary algebras associated to diagrams}
In this section we are going to associate to a given (oriented or unoriented) diagram a certain universal ternary algebra based on the conditions from the last two sections. We will prove the invariance of such algebras, up to isomorphism, with respect to the Reidemeister moves. We begin with two preliminary lemmas that can be found in \cite{Coh}; this book also contains the definitions of some elementary notions of universal algebra that appear here. First, we briefly recall the idea of a presentation of an algebra. 

Let $\Omega$ be an operator domain, and $\mathcal{A}$ be an algebra in a given category $\mathfrak{C}$ of $\Omega$-algebras. Instead of considering the carrier and the multiplication tables of $\mathcal{A}$, we could find a generating set $X$ of $\mathcal{A}$, and a family of relations $u=v$, where $u$ and $v$ are certain $\Omega$-words in $x_1,\ldots,x_n\in X$, that suffices to determine the effect of all the operators in $\mathcal{A}$. Then for a set $R$ of such relations, the definition of $\mathcal{A}$ using $X$ and $R$
is called a {\it presentation} of $\mathcal{A}$, and is denoted by $\mathfrak{C}\{X|R\}$. For the sake of convenience, this notion is extended: $X$ is taken to be a set of symbols such that each symbol is identified with an element of $\mathcal{A}$, and the corresponding elements of $\mathcal{A}$ generate $\mathcal{A}$. Thus, distinct elements of $X$ may represent the same element of $\mathcal{A}$.

Any variety $\mathcal{V}$ of $\Omega$-algebras gives a subcategory $\mathfrak{V}$ of the category $(\Omega)$ of all $\Omega$-algebras. It is simply the full subcategory of $(\Omega)$ whose objects are members of $\mathcal{V}$. Such $\mathfrak{V}$ is always a residual category, and the following two lemmas apply.

\begin{lemma}
Let 
$\mathfrak{C}$ be a residual category of $\Omega$-algebras. Then every presentation $\mathfrak{C}\{X|R\}$ defines a $\mathfrak{C}$-algebra, unique up to isomorphism.
\end{lemma}

A presentation $\mathfrak{C}\{X|R\}$ is called {\it finite} if both $X$ and $R$ are finite.

\begin{lemma}[Tietze]\label{Tietze}
Let 
$\mathfrak{C}$ be a residual category of $\Omega$-algebras, and let $\mathfrak{C}\{X|R\}$ be a finite presentation of a 
$\mathfrak{C}$-algebra 
$\mathcal{A}$. Then any other finite presentation of $\mathcal{A}$ is obtained from it by the following operations and their inverses:\\
(TI) If $(u,v)$ is a consequence of $R$, replace $R$ by $R\cup \{(u,v)\}$. Here we used the notation $(u,v)$ instead of $u=v$; a consequence can be described precisely as an element of the smallest congruence containing $R$.\\
(TII) If $u$ is any word in $X$ and $y$ is any letter not occurring in $X$, replace $X$ by $X\cup\{y\}$ and $R$ by $R\cup \{(y,u)\}$.
\end{lemma}

We will now proceed to defining algebras associated to link diagrams.
\begin{definition}
Let $\mathcal{V}^{ch}$ be the variety of algebras of a type (3,3) with operator symbols $B$ and $W$, satisfying the universal relations:
\begin{enumerate}
\item[(1)] $b(acbB)aB=c$ 
\item[(2)] $b(acbW)aW=c$                              
\item[(3)] $ba(abcB)W=c=ba(abcW)B$               
\item[(4)] $(cabB)baW=c=(cabW)baB$ 
\item[(5)] $(abcW)cdB=[ab(bcdB)W](bcdB)dB$
\item[(6)] $ab(bcdB)W=a(abcW)[(abcW)cdB]W$
\item[(7)] $(abcB)cdW=[ab(bcdW)B](bcdW)dW$
\item[(8)] $ab(bcdW)B=a(abcB)[(abcB)cdW]B$.
\end{enumerate}
Let $\mathfrak{V}^{ch}$ be a full subcategory of $(\Omega)$ with the objects $\mathcal{V}^{ch}$. Let $X$ be the set of symbols corresponding to the regions of a given unoriented diagram $D$ of a link $L$. Let the diagram be checkerboard-colored in such a way that the color of the outside region is white. The set of relations $R$ corresponds to crossings in the following way: at a crossing, choose one of the four regions, and let $q$ be its symbol; the relation assigned to the crossing is $q=xyzT$, where $x,y$, and $z$ are the symbols of the remaining regions near the crossing taken counter-clockwise if the region with label $q$ contains the positive marker, and clockwise otherwise; here $T=W$ if the region labeled $q$ is white, and $T=B$ if it is black. Thus, when looking at the left crossing on the Fig. \ref{unorientedop}, one of the relations $a=dcbB$, $b=cdaW$, $c=badB$ or $d=abcW$ needs to be included in $R$, the other three are its consequences assuming the axioms (1)-(4). For the crossing on the right of Fig. \ref{unorientedop}, the relation is chosen among:
$a=dcbW$, $b=cdaB$, $c=badW$, and $d=abcB$. If $R$ is the set of relations taken over all crossings, then we define the ternary algebra of the diagram $D$, denoted by $\mathcal{T}^{ch}(D)$, as the algebra with presentation $\mathfrak{V}^{ch}\{X|R\}$. Its isomorphism class is an invariant of the link, so we can write $\mathcal{T}^{ch}(L)$, and call it the {\it ternary algebra of the unoriented link $L$}.
\end{definition}

\begin{remark}
We have seen in the second section that taking $abcW=abcB=ab^{-1}c$ gives the relations from the fundamental group, and taking $abcB=c^{-1}ab$, $abcW=bac^{-1}$, produces the relations from the core group. Thus, the algebra $\mathcal{T}^{ch}(L)$ contains the information from these groups.
\end{remark}

We define similar algebras for oriented links.

\begin{definition}
Let $\mathcal{V}$ be the variety of algebras of a type (3,3) with operator symbols $C$ and $S$, satisfying the universal relations:
\begin{enumerate}
\item[(1)] $b(acbC)aC=c$ 
\item[(2)] $b(acbS)aS=c$                              
\item[(3)] $ba(abcC)S=ba(abcS)C=c$               
\item[(4)] $(cabC)baS=(cabS)baC=c$ 
\item[(5)] $(abcC)cdC=[ab(bcdC)C](bcdC)dC$
\item[(6)] $ab(bcdC)C=a(abcC)[(abcC)cdC]C$.
\end{enumerate}
Let $\mathfrak{V}$ be a full subcategory of $(\Omega)$ given by $\mathcal{V}$. Let $X$ be the set of symbols corresponding to the regions of a given oriented diagram $D$ of a link $L$. The set of relations $R$ corresponds to the crossings of $D$ as follows: at a positive crossing like the one on the left of 
Fig. \ref{orientedop}, choose one of the relations:
$a=dcbS$, $b=cdaC$, $c=badS$, or $d=abcC$. The remaining three relations will follow from axioms (1)-(4). For a negative crossing as on the right of Fig. \ref{orientedop}, select one relation among: $a=dcbC$, $b=cdaS$, $c=badC$, and $d=abcS$. In other words, the relations from Definition \ref{col_or} are used, one for each crossing. If $R$ is the set of all such relations, then we define the ternary algebra of the diagram $D$, denoted by $\mathcal{T}(D)$, as the algebra with presentation $\mathfrak{V}\{X|R\}$. As we prove below, its isomorphism class is an invariant of the link, so we can write $\mathcal{T}(L)$, and call it the {\it ternary algebra of the oriented link $L$}.
\end{definition}

\begin{theorem}
The isomorphism classes of algebras $\mathcal{T}^{ch}(L)$ and $\mathcal{T}(L)$ are invariants of links.
\end{theorem}
\begin{proof}
We will prove the theorem in the case of algebra $\mathcal{T}(L)$, the proof for $\mathcal{T}^{ch}(L)$ is similar. We use Theorem \ref{Tietze} to show that the presentations assigned to diagrams before and after the Reidemeister moves give the same algebra, up to isomorphism.

For the first Reidemeister move, depicted in Fig. \ref{R1or}, adding a new generator $c$ and expressing it with old generators via $c=babC$ is just an instance of the Tietze operation (TII). 

After the second Reidemeister move, two new relations: $(i)\ x=abcS$ and $(ii)\ y=cxaS$ are added to the presentation, see Fig. \ref{R2or}. Note that $x$ appears only in these two relations, and $y=b$ is their consequence, so we can add it to the set of relations using (TI), and then remove $(ii)$ with (TI)$^{-1}$ because it is now a consequence of $y=b$ and $(i)$. Then replace $y$ by $b$ in all the relations except $y=b$ (using (TI)'s and (TI)$^{-1}$'s). Finally, using the operations (TII)$^{-1}$, remove $y$ and then $x$, together with the relations $y=b$ and $(i)$, to obtain the presentation from before the move. 

Now we show the equivalence of presentations before and after the third Reidemeister move (Fig. \ref{R3or}).
On the left of Fig. \ref{R3or}, the associated presentation $\mathcal{T}(D)$ contains the relations: $(i)\ e=bcdC$, $(ii)\ f=abeC$, and $(iii)\ g=fedC$. First, add the consequences $(iv)\ f=ab(bcdC)C$ and $(v)\ g=[ab(bcdC)C](bcdC)dC$ using operations (TI), and remove $(ii)$ and $(iii)$ which are now consequences of $(i)$, $(iv)$, and $(v)$. Now, the generator $e$ and the relation $(i)$ can be removed with (TII)$^{-1}$. Next, replace $f$ by 
$ab(bcdC)C$, and $g$ by $[ab(bcdC)C](bcdC)dC$ in all the relations except $(iv)$ and $(v)$, and remove these two generators, together with $(iv)$ and $(v)$, using operations (TII)$^{-1}$. Call the resulting presentation $P_1$. We perform analogous operations on the presentation assigned to the diagram whose part is represented by the right side of the Fig. \ref{R3or}. The relations include: 
$(i')\ h=abcC$, $(ii')\ i=hcdC$, and $(iii')\ j=ahiC$. The consequences $(iv')\ i=(abcC)cdC$ and $(v')\ j=a(abcC)[(abcC)cdC]C$ are added, relations 
$(ii')$, $(iii')$, and then the generator $h$ and the relation $(i')$ are removed. In all the relations except $(iv')$ and $(v')$, $i$ gets replaced by $(abcC)cdC$,
and $j$ by $a(abcC)[(abcC)cdC]C$. Finally, remove $i$ and $j$, together with $(iv')$ and $(v')$, using operations (TII)$^{-1}$. Call this presentation $P_2$. The equivalence of $P_1$ and $P_2$ follows from the fact that we work in a variety of algebras satisfying axioms (5) and (6). 
\end{proof}

\begin{lemma}\label{colshom}
Any checkerboard ternary coloring of a diagram $D$ using the elements of an algebra $\mathcal{A}$ can be identified with a homomorphism from 
$\mathcal{T}^{ch}(D)$ to $\mathcal{A}$. Any ternary coloring of an oriented diagram $D$ using algebra $\mathcal{A}$ can be identified with a homomorphism from $\mathcal{T}(D)$ to $\mathcal{A}$.
\end{lemma}

The proof follows from the following lemma \cite{Coh}.

\begin{lemma}
Let $A$ and $B$ be any $\Omega$-algebras. Given a set $X$ and mappings $\alpha\colon X\to A$ and $\beta\colon X\to B$ such that:\\
(i) $A$ is generated by $im\ \alpha$,\\
(ii) any relation in $A$ between the elements $x\alpha$  $(x\in X)$ also holds between the corresponding elements $x\beta$ in $B$;\\
then there exists a unique homomorphism $\phi\colon A\to B$ such that $\alpha\phi=\beta$.
\end{lemma}

Lemma \ref{colshom} follows if we take $X$ to be the set of symbols assigned to the regions of the diagram $D$, $A=\mathcal{T}^{ch}(D)$
(or $A=\mathcal{T}(D)$ in the oriented case), and $B=\mathcal{A}$.
                          
\section{Search for ternary operators using nonassociative binary operations}
In this section we describe the results of a search for ternary operations involving binary operations from nonassociative structures. It turned out that
some varieties of loops contain interesting examples. Before listing the results, we recall preliminary definitions.

\begin{definition}
A {\it quasigroup} is a groupoid $(Q,*)$ such that the equation $x*y=z$ has a unique solution in $Q$ whenever two
of the three elements $x$, $y$, $z$ of $Q$ are specified. Equivalently, a quasigroup can be defined as a set $Q$ with three binary operations $*$, $\bs$, and $/$, satisfying the identities: $x\bs(x*y)=y$, $x*(x\bs y)=y$, $(x*y)/y=x$, and $(x/y)*y=x$, for any $x$, $y\in Q$. A {\it loop} $L$ is a quasigroup with an identity element $e$ such that $x*e=x=e*x$, for all $x\in L$. Standard references for the 
loop theory include \cite{Bru} and \cite{Pfl}.
\end{definition}

\begin{definition}
A {\it left Bol loop} is a loop $L$ satisfying the condition \[x*(y*(x*z))=(x*(y*x))*z,\] for all $x$, $y$, and $z\in L$.
\end{definition}

In a Bol loop $L$, the subloop generated by any element $x\in L$ is a group. We can therefore define $x^{-1}$ as the inverse of $x$ in that group.
Left Bol loops satisfy the {\it left inverse property}: $x^{-1}*(x*y)=y$, for any $x$ and $y\in L$; then we can write $x\bs y=x^{-1}y$. They are also 
{\it left alternative}: $x*(x*y) = (x*x)*y$; see \cite{Rob} for the material on Bol loops.

\begin{definition}
A {\it Moufang loop} is a loop $L$ satisfying one of the following equivalent identities:
\begin{enumerate}
\item[(1)] $z*(x*(z*y))=((z*x)*z)*y,$
\item[(2)] $x*(z*(y*z))=((x*z)*y)*z,$
\item[(3)] $(z*x)(y*z)=(z*(x*y))*z,$
\item[(4)] $(z*x)(y*z)=z*((x*y)*z),$
for all $x$, $y$, and $z\in L$.
\end{enumerate}
\end{definition}

Every Moufang loop is a left (and right) Bol loop, in particular it is {\it right alternative}: $x*(y*y)=(x*y)*y$. In addition to the left inverse property, they have the {\it right inverse property} $(y*x)*x^{-1}=y$,
and it follows that $y/x=y*x^{-1}$. These two inverse properties imply the {\it antiautomorphic inverse property}: $(x*y)^{-1}=y^{-1}*x^{-1}$. Moufang loops are always {\it flexible}: $x(yx)=(xy)x$, which is a consequence of the fact that they are {\it di-associative}, i.e., the subloop generated by any two elements is 
a group; because of that, the parenthesis in expressions involving only two elements are often dropped.

\begin{definition}
A loop $L$ is {\it conjugacy closed} if it satisfies the equations:
\begin{enumerate}
\item[(1)] $(x*y)*z=(x*z)*(z\bs(y*z)),$
\item[(2)] $z*(y*x)=((z*y)/z)*(z*x),$ for all $x$, $y$, and $z\in L$. 
\end{enumerate}
\end{definition}
For information about conjugacy closed loops see, for example, \cite{KK2}.

\begin{definition}
An {\it extra loop} is a loop $L$ satisfying one of the following equivalent conditions:
\begin{enumerate}
\item[(1)] $(x*(y*z))*y=(x*y)*(z*y),$
\item[(2)] $(y*z)*(y*x)=y*((z*y)*x),$
\item[(3)] $((x*y)*z)*x=x*(y*(z*x)),$
for all $x$, $y$, and $z\in L$. 
\end{enumerate}
\end{definition}
Every extra loop is a Moufang loop, and is conjugacy closed. It is also a {\it C-loop}, i.e., it satisfies $x*(y*(y*z))=((x*y)*y)*z$, which implies two conditions:
\begin{enumerate}
\item[(LC)] $(x*x)*(y*z)=(x*(x*y))*z,$
\item[(RC)] $x*((y*z)*z)=(x*y)*(z*z).$
\end{enumerate}
Note that any group is an extra loop. For more information about extra loops see \cite{KK}.
The following is an example of an extra loop that is not a group (see \cite{PV} and the references therein): let $G$ be a group, and $M(G,2)$ be the set 
$G\times\{0,1\}$ equipped with the operation
$(g,0)*(h,0)=(gh,0)$, $(g,0)*(h,1)=(hg,1)$, $(g,1)*(h,0)=(gh^{-1},1)$, and $(g,1)*(h,1)=(h^{-1}g,0)$. Then $M(G,2)$ is a nonassociative Moufang loop if and only if $G$ is nonabelian, and $M(D_4,2)$, where $D_4$ is a dihedral group with eight elements, is an extra loop.

The lists presented below contain pairs of ternary operators that can be taken as $C$ and $S$ (or $S$ and $C$) in the oriented case, and as $W$ and $B$ (or $B$ and $W$) in the unoriented case. They come from Moufang loops, and extra loops. We searched for the candidates using GAP, and in our program we used the libraries of loops, and some functions from the GAP package \cite{loops}. Then the candidates were verified by hand to satisfy the oriented case axioms (1)-(6), or unoriented case axioms (1)-(8). The operators for the oriented case turned out to be always the same in each pair, so these pairs do not distinguish the orientation of diagrams; perhaps it's better to look for invariants distinguishing orientation among ternary operators that do not come from binary operations. We note, however, that since loops generalize groups, some of the listed operators give generalizations of the Dehn presentation of the knot group to the category of extra loops, and some give generalizations of the core group presentations that we considered in Section \ref{unoriented}. This last fact is not surprising, as it is well known that the cores of Moufang loops are involutory quandles \cite{Bru}.

\noindent Moufang loops - oriented case:
\begin{enumerate}
\item[(m1)] $(b*a^{-1})*c,\ (b*a^{-1})*c$
\item[(m2)] $(b*c^{-1})*a,\ (b*c^{-1})*a$
\item[(m3)] $a*(c^{-1}*b),\ a*(c^{-1}*b)$
\item[(m4)] $c*(a^{-1}*b),\ c*(a^{-1}*b)$
\end{enumerate}
\noindent Moufang loops - unoriented case, in addition to the pairs (m1)-(m4), we have:
\begin{enumerate}
\item[(m5)] $a^{-1}*(c*b),\ (b*c)*a^{-1}$
\item[(m6)] $c^{-1}*(a*b),\ (b*a)*c^{-1}$
\end{enumerate}

By comparing the above pairs (m1)-(m6) with the pairs (g1)-(g9) from section \ref{unoriented}, we can notice that on the level of groups they would give core group relations. As an example, we will prove that (m1) satisfies the required conditions.

\begin{lemma}
The pair $xyzS=xyzC=(y*x^{-1})*z$ satisfies the oriented case identities (1)-(6) in the category of Moufang loops.
\end{lemma}
\begin{proof}
Because the operators are equal, we need to prove only conditions (1), (3), (4), (5), and (6). We will suppress `$*$' in the calculations, but we will use dots to separate the inputs of operators.
\begin{align*}
(1)\ &b.(a.c.bC).aC=b.(ca^{-1})b.aC=(((ca^{-1})b)b^{-1})a=c;\\
(3)\ &b.a.(a.b.cC)C=b.a.(ba^{-1})cC=(ab^{-1})((ba^{-1})c)=
     (ba^{-1})^{-1}((ba^{-1})c)=c;\\
(4)\ &(c.a.bC).b.aC=(ac^{-1})b.b.aC=(b((ac^{-1})b)^{-1})a=(b(b^{-1}(ac^{-1})^{-1}))a=\\
     &(ca^{-1})a=c;\\
(5)\ &(a.b.cC).c.dC=(ba^{-1})c.c.dC=(c((ba^{-1})c)^{-1})d=(c(c^{-1}(ba^{-1})^{-1}))d=\\
     &(ab^{-1})d,\\
     &[a.b.(b.c.dC)C].(b.c.dC).dC=(a.b.(cb^{-1})dC).(cb^{-1})d.dC=\\
     &(ba^{-1})((cb^{-1})d).(cb^{-1})d.dC=(((cb^{-1})d)((ba^{-1})((cb^{-1})d))^{-1})d=(ab^{-1})d;
\end{align*}
In the proof of identity (6), we use the substitutions: $k=ba^{-1}$, and $q=a^{-1}k^{-1}$.
\begin{align*}
    &a.b.(b.c.dC)C=a.b.(cb^{-1})dC=(ba^{-1})((cb^{-1})d)=k((c(a^{-1}k^{-1}))d)=k((cq)d),\\
    &a.(a.b.cC)[(a.b.cC).c.dC]C=a.(ba^{-1})c.[(ba^{-1})c.c.dC]C=\\
    &a.(ba^{-1})c.(c(c^{-1}(ab^{-1})))dC=a.(ba^{-1})c.(ab^{-1})dC=\\
    &(((ba^{-1})c)a^{-1})((ab^{-1})d)=((kc)a^{-1})(k^{-1}d)=\\
    &((kc)(qk))(k^{-1}d)=(k(cq)k)(k^{-1}d)=k((cq)(k(k^{-1}d)))=k((cq)d).
\end{align*}
\end{proof}

\noindent Extra loops - oriented case:
\begin{enumerate}
\item[(e1)] $(a*b^{-1})*c,\ (a*b^{-1})*c$
\item[(e2)] $(a*c^{-1})*b,\ (a*c^{-1})*b$
\item[(e3)] $(b*a^{-1})*c,\ (b*a^{-1})*c$
\item[(e4)] $(b*c^{-1})*a,\ (b*c^{-1})*a$
\item[(e5)] $(c*a^{-1})*b,\ (c*a^{-1})*b$
\item[(e6)] $(c*b^{-1})*a,\ (c*b^{-1})*a$
\item[(e7)] $a*(b^{-1}*c),\ a*(b^{-1}*c)$
\item[(e8)] $a*(c^{-1}*b),\ a*(c^{-1}*b)$
\item[(e9)] $b*(a^{-1}*c),\ b*(a^{-1}*c)$
\item[(e10)] $b*(c^{-1}*a),\ b*(c^{-1}*a)$
\item[(e11)] $c*(a^{-1}*b),\ c*(a^{-1}*b)$
\item[(e12)] $c*(b^{-1}*a),\ c*(b^{-1}*a)$
\end{enumerate}

\noindent Extra loops - unoriented case, in addition to the pairs (e1)-(e12), we have:
\begin{enumerate}
\item[(e13)] $(a^{-1}*c)*b,\ b*(c*a^{-1})$
\item[(e14)] $(c^{-1}*a)*b,\ b*(a*c^{-1})$
\item[(e15)] $a^{-1}*(c*b),\ (b*c)*a^{-1}$
\item[(e16)] $c^{-1}*(a*b),\ (b*a)*c^{-1}$
\item[(e17)] $(a^{-1}*b^{-1})*c^{-1},\ c^{-1}*(b^{-1}*a^{-1})$
\item[(e18)] $a^{-1}*(b^{-1}*c^{-1}),\ (c^{-1}*b^{-1})*a^{-1}$
\end{enumerate}

We will now prove that the pair (e1), which on the level of groups yields the Dehn presentation of the knot group, satisfies the required operator conditions. 

\begin{lemma}
The pair $xyzS=xyzC=(xy^{-1})z$ satisfies the oriented case identities (1)-(6) in the category of extra loops.
\end{lemma}
\begin{proof}
\begin{align*}
(1)\ &b.(a.c.bC).aC=b.(ac^{-1})b.aC=(b((ac^{-1})b)^{-1})a=(b(b^{-1}(ac^{-1})^{-1}))a=\\
     &(ca^{-1})a=c;\\
(3)\ &b.a.(a.b.cC)C=b.a.(ab^{-1})cC=(ba^{-1})((ab^{-1})c)=(ab^{-1})^{-1}((ab^{-1})c)=c;\\
(4)\ &(c.a.bC).b.aC=(ca^{-1})b.b.aC=(((ca^{-1})b)b^{-1})a=c;\\
(5)\ &(a.b.cC).c.dC=(ab^{-1})c.c.dC=(((ab^{-1})c)c^{-1})d=(ab^{-1})d,\\
     &[a.b.(b.c.dC)C].(b.c.dC).dC=(a.b.(bc^{-1})dC).(bc^{-1})d.dC=\\
     &(ab^{-1})((bc^{-1})d).(bc^{-1})d.dC=(((ab^{-1})((bc^{-1})d))((bc^{-1})d)^{-1})d=(ab^{-1})d;
\end{align*}
To prove identity (6), first we will prove that $(z(bc^{-1}))z^{-1}=(zb)(c^{-1}z^{-1})$, for any elements $b$, $c$, $z$.
From the condition (RC) we have: \[(zb)(c^{-1}z^{-1})=(zb)(((c^{-1}z)z^{-1})z^{-1})=((zb)(c^{-1}z))(z^{-1}z^{-1}).\] From the Moufang identity (3):
\[((zb)(c^{-1}z))(z^{-1}z^{-1})=((z(bc^{-1}))z)(z^{-1}z^{-1}).\] Finally, from the right alternative property, we have:
\[((z(bc^{-1}))z)(z^{-1}z^{-1})=(((z(bc^{-1}))z)z^{-1})z^{-1}=(z(bc^{-1}))z^{-1}.\]
The left side of (6) gives:
\[ a.b.(b.c.dC)C=ab((bc^{-1})d)C=(ab^{-1})((bc^{-1})d),\]
and on the right side, we have:
\begin{align*}
&a.(a.b.cC)[(a.b.cC).c.dC]C=a.(ab^{-1})c.[(ab^{-1})c.c.dC]C=\\
&a.(ab^{-1})c.(((ab^{-1})c)c^{-1})dC=a.(ab^{-1})c.(ab^{-1})dC=\\
&(a((ab^{-1})c)^{-1})((ab^{-1})d)=(a(c^{-1}(ab^{-1})^{-1}))((ab^{-1})d).
\end{align*}
Now the equality of two sides follows from the conjugacy closed condition (2): $z(yx)=((zy)z^{-1})(zx),$ where we take $z=ab^{-1}$, $x=d$, $y=bc^{-1}$, and from the already proven $(z(bc^{-1}))z^{-1}=(zb)(c^{-1}z^{-1})$:
\begin{align*}
&(ab^{-1})((bc^{-1})d)=z(yx)=((zy)z^{-1})(zx)=((z(bc^{-1}))z^{-1})(zd)=\\
&((zb)(c^{-1}z^{-1}))(zd)=(a(c^{-1}(ab^{-1})^{-1}))((ab^{-1})d).
\end{align*}

\end{proof}

\noindent Computer computations suggest that the following pairs work in the case of left Bol loops.\\
\noindent Left Bol loops - oriented case:
\begin{enumerate}
\item[(b1)] $(b/a)*c,\ (b/a)*c$
\item[(b2)] $(b/c)*a,\ (b/c)*a$
\end{enumerate}
\noindent Left Bol loops - unoriented case, in addition to the pairs (b1) and (b2), there are pairs:
\begin{enumerate}
\item[(b3)] $(b/a^{-1})*c^{-1},\ ((a/b^{-1})\bs c)^{-1}$
\item[(b4)] $(b/c^{-1})*a^{-1},\ ((c/b^{-1})\bs a)^{-1}$
\end{enumerate}

We also note that there are more binary structures with a multitude of operators satisfying only the `distributivity' conditions (5)-(8) in the unoriented case, and (5)-(6) in the oriented case.

\section{Concluding remarks}
The purpose of this paper was to show that ternary operations can be explicitly used in knot theory in a simple way, and that they arise naturally from known invariants. This connection was the reason why we worked with the four regions around the crossing rather than with the four semi-arcs of the crossing. We believe that the full strength of ternary operations will become apparent in the case of knotted surfaces, and, more generally, $n$-ary operations could be successfully used for the higher dimensional knots. Naturally, the homology theory for the corresponding $n$-ary structures could prove to be of great use.

\end{document}